\documentclass[12pt]{article}

\usepackage[table]{xcolor}
\usepackage{amsmath, amsthm}
\usepackage{amsfonts,amssymb}
\usepackage{hyperref}
\usepackage{authblk}
\usepackage[all]{xy}
\usepackage{enumitem}
\usepackage{stmaryrd}
\usepackage[makeroom]{cancel}
\usepackage{bbm}
\usepackage{fancyvrb}
\usepackage[linesnumbered,lined,boxed]{algorithm2e}
\usepackage{xstring}
\usepackage{tikz}
\usetikzlibrary{intersections}

\newcolumntype{C}{>{\hfil$}p{7.5mm}<{$\hfil}}
\newenvironment{Smallmatrix}[1]
  {\arraycolsep=3pt\scriptsize
   \array{CCC}{#1}}
  {\endarray }

\renewcommand{\le}{\leqslant}
\renewcommand{\ge}{\geqslant}
\newcommand{\N}{\mathbb{N}}
\newcommand{\Z}{\mathbb{Z}}
\newcommand{\Q}{\mathbb{Q}}

\newcommand{\F}{\mathbb{F}}

\renewcommand{\P}{\mathbb{P}}

\newcommand{\A}{\mathbb{A}}

\newcommand{\Gal}{\operatorname{Gal}}

\newcommand{\GL}{\operatorname{GL}}
\newcommand{\SL}{\operatorname{SL}}

\newcommand{\PGL}{\operatorname{PGL}}

\newcommand{\U}{\operatorname{U}}
\newcommand{\SU}{\operatorname{SU}}
\newcommand{\PSU}{\operatorname{PSU}}
\newcommand{\PSUt}{\PSU_3(\F_9)}
\newcommand{\SUt}{\SU_3(\F_9)}
\newcommand{\Ut}{\U_3(\F_9)}
\newcommand{\Fl}{{\F_\ell}}

\newcommand{\Ker}{\operatorname{Ker}}

\newcommand{\Frob}{\operatorname{Frob}}

\newcommand{\Jac}{\operatorname{Jac}}
\renewcommand{\H}{\operatorname{H}_\text{\'et}}

\newcommand{\smat}[4]{\left[ \begin{smallmatrix} #1 & #2 \\ #3 & #4 \end{smallmatrix} \right]}

\numberwithin{equation}{section}

\newtheorem{thm}[equation]{Theorem}

\theoremstyle{definition}

\newtheorem{rk}[equation]{Remark}

\makeatletter
\newcommand{\subjclass}[2][2010]{%
  \let\@oldtitle\@title%
  \gdef\@title{\@oldtitle\footnotetext{#1 \emph{Mathematics subject classification:} #2}}%
}
\newcommand{\keywords}[1]{%
  \let\@@oldtitle\@title%
  \gdef\@title{\@@oldtitle\footnotetext{\emph{Key words and phrases.} #1.}}%
}
\let\c@table\c@equation
\let\c@figure\c@equation
\makeatother
\makeatother

\title{Explicit computation of a Galois representation \\ attached to an eigenform over~$\SL_3$ \\ from the~$\H^2$ of a surface}
\author{Nicolas Mascot\thanks{\href{mailto:nm116@aub.edu.lb}{nm116@aub.edu.lb}}}
\affil{\scriptsize{AUB, Beirut, Lebanon}}

\VerbatimFootnotes

\begin{document}

\maketitle

\begin{abstract}
We sketch a method to compute mod~$\ell$ Galois representations contained in the~$\H^2$ of surfaces. We apply this method to the case of a representation with values in~$\GL_3(\F_9)$ attached to an eigenform over a congruence subgroup of~$\SL_3$. We obtain in particular a polynomial with Galois group isomorphic to the simple group~$\PSUt$ and ramified at~2 and~3 only. 
\end{abstract}

\renewcommand{\abstractname}{Acknowledgements}
\begin{abstract}

The author thanks Jean-Marc Couveignes for his suggestion to apply the \emph{d\'evissage} principle to the case of the~$\H^2$ of surfaces, and in particular on the example presented in this article; Fran\c{c}ois Brunault for his help with the details of the proof of theorem~\ref{thm:devissage}; and Bert van Geemen for pointing out the possibility to read the monodromy around the bad fibres of an elliptic surface off their Kodaira types.

The computer algebra packages used for the computations presented in this article were~\cite{gp} and~\cite{Magma}. The computations were carried out on the Warwick mathematics institute computer cluster provided by the EPSRC Programme Grant EP/K034383/1
``LMF: L-Functions and Modular Forms''.
\end{abstract}


\newpage

\section{Introduction}\label{sect:intro}

Several techniques (such as~\cite{CE11},~\cite{companion}, and~\cite{Hensel}) have recently been developed to compute explicitly mod~$\ell$ Galois representations afforded in the torsion of Jacobians of curves. However, many ``interesting'' representations (e.g. in view of the Langlands program) are not naturally found in Jacobians, but rather in higher \'etale cohomology spaces of higher-dimensional varieties. They are thus inaccessible to the aforementioned methods, and, to the author's knowledge, computational methods to deal with these representations have not yet been developed.

The purpose of this article is to sketch such a method. Although it is still at an experimental stage, it is already sufficiently advanced for us to be able to prove our concept by giving a concrete example of application, which we also present in this article.

\begin{rk}
Very general algorithms to compute with \'etale cohomology are presented \linebreak in~\cite{Madore} and~\cite{Testa}; however, as far as we know these algorithms are not really practical and have never been implemented. Our purpose is to present a method which is really practical on a moderately simple case, and ought to be generalizable to other similar cases.
\end{rk}

\bigskip

The concrete example that we have chosen comes from~\cite{Top}, where B. van Geemen and J. Top give evidence towards a conjecture of L. Clozel's . They exhibit a Hecke eigenform~$u$ over a congruence subgroup of~$\SL_3(\Z)$ of level 128,
whose Hecke eigenvalues lie is~$\Z[2\sqrt{-1}]$, and an algebraic surface~$S$ over~$\Q$, such that for all primes~$\ell \in \N$, the~$\ell$-adic~$\operatorname{H}^2$ of~$S$ contains a Galois-submodule affording the quadratic twist by~$-2$ of the~$\ell$-adic representation
\[ \widetilde \rho_{u,\ell} : \Gal(\overline \Q / \Q) \longrightarrow \GL_3\big(\Q_\ell(\sqrt{-1})\big) \]
attached to~$u$.

\begin{rk}At the time when~\cite{Top} was published, representations attached to this kind of modular forms were not even known to exist, especially when, like~$\widetilde \rho_{u,\ell}$, they are not self-dual. The existence of these representations was established recently independently by~\cite{Harris} and~\cite{Scholze}, and the fact that~$\widetilde \rho_{u,\ell}$ is indeed afforded by the~$\ell$-adic~$\operatorname{H}^2$ of~$S$ is proved in~\cite{Ito}.
\end{rk}

As~$u$ has level~$128=2^7$, the~$\ell$-adic representation~$\widetilde \rho_{u,\ell}$ is unramified away from 2 and~$\ell$. According to~\cite{Top}, for each unramified~$p$, the characteristic polynomial of~$\widetilde \rho_{u,\ell}(\Frob_p)$ is
\begin{equation} \chi_p(x) = x^3-a_p x^2 +p \overline{a_p} x - p^3 \in \Q_\ell(\sqrt{-1})[x], \label{eqn:chip} \end{equation}
where~$a_p$ is the corresponding Hecke eigenvalue of~$u$, and~$\overline{a_p}$ is the image of~$a_p$ under complex conjugation. In particular, the determinant of this representation is  the cube of the~$\ell$-adic cyclotomic character, and the value of~$a_p$ can be recovered as the trace of the Frobenius.

\bigskip

\begin{rk}\label{rk:rho2triv}
Since~$a_p$ lies in~$\Z[2\sqrt{-1}]$ for all~$p$, the characteristic polynomial~$\chi_p(x)$ is always congruent to~$(x-1)^3$ mod~$2$. This shows that the mod~$2$ representation is trivial (up to semi-simplification). Therefore, we have chosen to consider the more interesting (and challenging) case~$\ell=3$.
\end{rk}

\pagebreak

The purpose of this article is thus to explain how we have almost certainly succeeded to computed explicitly the mod~$3$ representation
\[ \rho_{u,3} : \Gal(\overline \Q / \Q) \longrightarrow \GL_3(\F_9) \]
found up to twist in~$\H^2(S_{\overline \Q}, \Z/3 \Z)$.

\bigskip

Here and in the rest of the article, by \emph{computing explicitly a mod~$\ell$ Galois representation}, we mean computing a polynomial whose roots are permuted by Galois in the same way as the vectors in the space of the representation (so that its splitting field agrees with the field cut out by the representation), as well as extra data making it possible to determine for any unramified prime~$p \in \N$ the image of the the Frobenius at~$p$ up to conjugacy. In our case, this implies in particular that our computations allow us to determine mod 3 the eigenvalue~$a_p$ of~$u$ for all~$p \ge 5$. For instance, we can compute in four  seconds that if~$p=10^{1000}+453$ is the first prime after~$10^{1000}$, then~$a_p \equiv -1 \bmod 3 \Z[\sqrt{-1}]$. As a bonus, the Galois group of the polynomial that we thus obtain, which is by construction the image of~$\rho_{u,3}$, turns out to be a particularly interesting subgroup of~$\GL_3(\F_9)$ (cf. section~\ref{sect:image}).

\begin{rk}
Unfortunately, because of the reason given in remark~\ref{rk:bestappr}, we are unable to certify rigorously that the results of our computations are correct. However, the fact that we are eventually able to recover the value of the~$a_p$ mod 3 from the representation (cf. section~\ref{sect:Dok}) shows that our results are correct beyond reasonable doubt.
\end{rk}

The central idea making this computation possible, which we owe to J.-M. Couveignes, is a method to construct by \emph{d\'evissage} a curve~$C$ (depending on~$\ell$) such that the reduction mod~$\ell$ of~$\widetilde \rho_{u,\ell}$ is contained in the~$\ell$-torsion of the Jacobian of~$C$. The point is that once we have obtained an explicit model for~$C$, we are able (at least in theory) to compute the representation, thanks to our technique presented in~\cite{Hensel}. In principle, this \emph{d\'evissage} technique could be iterated to construct a curve whose Jacobian contains a given representation found in the~$\H^d$ of a variety of dimension~$d$.

\bigskip

We will sketch this construction in section~\ref{sect:devissage}, after which we will apply it to~$\rho_{u,3}$ in section~\ref{sect:curvered}. Next, we will explain in section~\ref{sect:Hensel} how we used the curve thus obtained to compute a polynomial corresponding to~$\rho_{u,3}$, after what we use this polynomial to determine the image of~$\rho_{u,3}$ in section~\ref{sect:image}. Finally, we will show in section~\ref{sect:Dok} how to compute the image of Frobenius elements.

\section{D\'evissage}\label{sect:devissage}

Suppose we are given a surface~$S$ defined over~$\Q$ as well as a prime~$\ell \in \N$ such \linebreak that~$\H^2(S_{\overline \Q},\Z/\ell\Z)$ contains a Galois-submodule affording a mod~$\ell$ Galois representation~$\rho$ that we wish to compute. We are going to show how to construct a curve~$C$ whose Jacobian will also contain~$\rho$ (up to twist) in its~$\ell$-torsion. Of course, this curve~$C$ will depend on~$\ell$. This construction is an example of the \emph{d\'evissage} method summarized in~\cite[3.4]{Deligne}.

\bigskip

Let~$\mu_\ell$ be the Galois module formed by the~$\ell$-th roots of unity. Given a Galois module~$M$ and an integer~$n \in \Z$, we will denote by~$M(n)$ the twist of~$M$ by the~$n$-th power of the mod~$\ell$ cyclotomic character. Thus~$\mu_\ell = (\Z/\ell\Z)(1)$ for instance. We will sometimes write~$\mu_\ell^\vee$ for~$(\Z/\ell\Z)(-1)$. Finally, we will also denote by~$\mu_\ell$ and~$(\Z/\ell\Z)(n)$ the corresponding constant sheaves on the \'etale site of a variety.


Recall~\cite[14.2]{Milne} that when~$X$ is a complete, connected and non-singular curve over~$\overline \Q$, we have canonical (and hence Galois-equivariant) identifications
\[ \H^r(X,\mu_\ell) \simeq \left\{ \begin{array}{ll} \mu_\ell & \text{ if } r=0, \\ J[\ell] & \text{ if } r=1, \\ \Z/\ell\Z & \text{ if } r=2,\end{array} \right. \]
where~$J$ is the Jacobian of~$X$. By tensoring out~$\mu_\ell$, we deduce the identifications
\begin{equation} \H^r(X,\Z/\ell\Z) \simeq \left\{ \begin{array}{ll} \Z/\ell\Z & \text{ if } r=0, \\ J[\ell](-1) & \text{ if } r=1, \\ \mu_\ell^\vee & \text{ if } r=2. \end{array} \right. \label{H_complete} \end{equation}

If we let~$U$ be~$X$ with finitely many points deleted, which is thus still a non-singular connected curved but is no longer complete, then we obtain
\begin{equation} \H^r(U,\Z/\ell\Z) \simeq \left\{ \begin{array}{ll} \Z/\ell\Z, & \text{ if } r=0, \\ J[\ell](-1) \text{ extended by copies of } \mu_\ell^\vee & \text{ if } r=1, \\ 0 & \text{ if } r=2, \end{array} \right. \label{H_incomplete} \end{equation}
where the first case is obvious, and the last two follow respectively from corollary~16.2 and proposition 14.12 of~\cite{Milne}.
\bigskip
%
%
 
Suppose now that we have a proper regular surface~$S$ over~$\Q$, equipped with a proper dominant morphism~$\pi : S \longrightarrow B$ to a non-singular complete curve~$B$ also defined over~$\Q$.

Let~$Z \subset B$ be a nonempty\footnote{We insist that~$Z$ must not be empty, because we will need~$B'$ to be affine later.} finite subset containing the image of the bad fibres of~$\pi$, and let~$Y = \pi^{-1}(Z) \subset S$. Define~$B'=B\setminus Z$, and~$S' = S \setminus Y$, so that the fibre~$S_b = S' \times_{B'} b$ at any~$b \in B'$ of the induced map~$\pi : S' \longrightarrow B'$ is a smooth proper curve. The representability of the relative Picard functor~\cite[9.4.4]{Neron} thus guarantees the existence of a cover~$\psi : C' \longrightarrow B'$ whose fibre at~$b \in B'$ is~$C'_b = \Jac(S_b)[\ell]$. 

The closed subscheme~$Y$ of~$S$ is made up of curves, possibly with multiplicities, and intersecting in some way. Define~$Y'$ to be the scheme obtained from~$Y$ by first passing to the reduced scheme structure, and then deleting the singular points. Thus$~Y'$ is a disjoint union of smooth curves. These curves need not be defined over $\Q$, and are thus permuted by Galois; let

\vspace{-4mm}

\[ \eta = \prod_{\text{Components of }Y'} \Fl \]
\vspace{-3mm}

\noindent be the corresponding mod $\ell$ permutation representation, and denote by $\eta(-1) = \eta \otimes \mu_\ell^\vee$ its twist by the inverse of mod $\ell$ cyclotomic character.

With this notation, we can prove that the ``interesting'' Galois representations which lie in~$\H^2(S_{\overline \Q},\Z/\ell\Z)$ are also afforded in~$\H^1(C'_{\overline \Q},\Z/\ell\Z)$:
\begin{thm}\label{thm:devissage}
Suppose~$\rho$ is a mod~$\ell$ Galois representation contained in~$\H^2(S_{\overline \Q},\Z/\ell\Z)$ (up to semi-simplification). Assume that $\rho$ has no Jordan-Hölder components of the form~$(\Z/\ell\Z)(n)$ for any~$n \in \Z$, and no component in common with $\eta(-1)$. Then the twist of~$\rho$ by the mod~$\ell$ cyclotomic character is also contained (up to semi-simplification) in~$\H^1(C'_{\overline \Q},\Z/\ell \Z)$.
\end{thm}

\begin{rk}
The number field cut out by $\eta(-1)$ is contained in the compositum of the $\ell$-th cyclotomic field and of the fields of definition of the components of the bad fibres of $\pi$. In general, we expect this field to be considerably smaller than that cut out by $\rho$ if $\rho$ is an ``interesting'' representation. For instance, for the surface considered in section \ref{sect:curvered} below, the field cut out by~$\eta$ is merely~$\Q(\sqrt2)$ (cf. remark \ref{rk:Kodaira}). Therefore, the requirement that $\rho$ have no common component with~$\eta(-1)$ ought to be harmless for ``interesting'' representations $\rho$.
\end{rk}

\begin{proof}


Let us first show that~$\rho$ is also contained in~$\H^2(S'_{\overline \Q}, \Z/\ell \Z)$, so that the bad fibres of~$\pi$ will no longer be a nuisance. The localization exact sequence~\cite[9.4]{Milne} shows that the kernel \linebreak of~$\H^2(S_{\overline \Q},\Z/\ell\Z) \longrightarrow \H^2(S'_{\overline \Q},\Z/\ell\Z)$ is a quotient of~$\operatorname{H}^2_{\text{\'et},Y}(S_{\overline \Q},\Z/\ell\Z)$. \'Etale cohomology with coefficients in~$(\Z/\ell\Z)(n)$ satisfies the Bloch-Ogus axioms~\cite{BlochOgus}, so in particular the Poincar\'e duality axiom~\cite[6.1.j]{Jannsen} shows that
\[ \operatorname{H}^2_{\text{\'et},Y}(S_{\overline \Q},\Z/\ell\Z) \simeq \operatorname{H}_2\!\big(Y_{\overline \Q},(\Z/\ell\Z)(2)\big). \]
 Applying~\cite[6.1.f]{Jannsen} twice shows that~$\operatorname{H}_2\!\big(Y_{\overline \Q},(\Z/\ell\Z)(2)\big) \simeq \operatorname{H}_2\!\big(Y'_{\overline \Q},(\Z/\ell\Z)(2)\big)$, and since~$Y'$ is a disjoint union of smooth curves, applying Poincar\'e duality~\cite[6.1.j]{Jannsen} and then~\eqref{H_complete} or~\eqref{H_incomplete} component-wise reveals that
 \[ \operatorname{H}_2\!\big(Y'_{\overline \Q},(\Z/\ell\Z)(2)\big) \simeq \eta(-1). \]
 Our assumptions on~$\rho$ thus show that it must be contained in~$\H^2(S'_{\overline \Q}, \Z/\ell \Z)$ as claimed.


\bigskip

Consider now the Leray spectral sequence~\cite[12.7]{Milne}
\[ E_2^{p,q} = \H^p(B'_{\overline \Q}, R^q \pi_* \Z/\ell \Z)  \Rightarrow \H^{p+q} (S'_{\overline \Q}, \Z/\ell \Z) \]
attached to~$\pi : S' \longrightarrow B'$. We know by proper base change~\cite[17.7]{Milne} that \[ R^q \pi_* \Z/\ell \Z = \H^q(S_b, \Z/\ell \Z), \] where by abuse of notation we denote by~$\mathcal{M}_b$ instead of~$\mathcal{M}$ the sheaf on~$B'$ whose stalk at~$b$ is~$\mathcal{M}_b$. Besides, the base~$B'$ and the fibres~$S_b$ are non-singular connected curves, so~$E_2^{p,q} = 0$ unless~$0 \le p,q \le 2$. Therefore~$E_2^{p,q} = E_\infty^{p,q}$ for all~$p,q$ such that~$p+q=1$, which means that~$H^2(S_{\overline \Q}, \Z/\ell\Z)$ admits a filtration with components
\begin{itemize}
\item~$\H^2(B'_{\overline \Q}, \Z/\ell\Z) =  0$ by~\eqref{H_incomplete},
\item~$\H^0\!\big(B'_{\overline \Q}, H^2(S_b,\Z/\ell\Z)\big) = \H^0(B'_{\overline \Q}, \mu_\ell^\vee) = \mu_\ell^\vee$ by~\eqref{H_complete} and~\eqref{H_incomplete},
\item and~$\H^1(B'_{\overline \Q}, \mathcal{F})$,
\end{itemize}
where~$\mathcal{F}$ is the sheaf on~$B'$ with stalks
\[ \mathcal{F}_b = \H^1(S_b, \Z/\ell\Z) = C'_b(-1) \]
by~\eqref{H_complete} and the definition~$C'_b = \Jac(S_b)[\ell]$. Our assumptions on~$\rho$ thus show that it must be contained in~$\H^1(B'_{\overline \Q}, \mathcal{F})$.

Similarly, the Leray spectral sequence 
\[ \H^p(B'_{\overline \Q}, R^q \psi_* \Z/\ell \Z)  \Rightarrow \H^{p+q} (C'_{\overline \Q}, \Z/\ell \Z) \]
attached to~$\psi : C' \longrightarrow B'$ shows that~$\H^1(C'_{\overline \Q}, \Z/\ell\Z) = \H^1(B'_{\overline \Q}, \mathcal{G})$, where
\[ \mathcal{G}_b = \H^0(C'_b, \psi_* \Z/\ell\Z) = \Fl^{C'_b}. \]
Since~$C'_b = \Jac(S_b)[\ell]$ is an Abelian~$\ell$-torsion group, there is a natural surjection
\[ \begin{array}{rcl} \mathcal{G}_b = \Fl^{C'_b} & \longrightarrow & C'_b = \mathcal{F}_b(1) \\ \lambda & \longmapsto & \sum_{c \in C'_b} \lambda_c \, c \end{array} \]
which, after twisting, yields an epimorphism~$\xymatrix{\mathcal{G}(-1) \ar@{->>}[r] & \mathcal{F}}$. Let~$\mathcal{K}$ by its kernel, so that we have the short exact sequence
\[ 0 \longrightarrow \mathcal{K} \longrightarrow \mathcal{G}(-1) \longrightarrow \mathcal{F} \longrightarrow 0 \]
of sheaves on~$B'_{\overline \Q}$. The associated long exact sequence contains
\[ \cdots \longrightarrow \H^1\!\big(B'_{\overline \Q}, \mathcal{G}(-1) \big) \longrightarrow \H^1(B'_{\overline \Q}, \mathcal{F}) \longrightarrow \H^2(B'_{\overline \Q}, \mathcal{K}) \longrightarrow \cdots, \]
and as~$\H^2(B'_{\overline \Q}, \mathcal{K}) = 0$ by~\cite[1.3.3.6.ii]{Deligne}, we conclude that since~$\rho$ is contained in~$\H^1(B'_{\overline \Q}, \mathcal{F})$, it also appears in 
\[ \H^1\!\big(B'_{\overline \Q}, \mathcal{G}(-1) \big)= \H^1(C'_{\overline \Q},\Z/\ell \Z)(-1). \]

%

\end{proof}

The curve~$C'$ constructed in theorem~\ref{thm:devissage} will not, in general, be connected, because of the zero section~$0 \in \Jac(S_b)[\ell] = C'_b$. However, we can modify the definition of~$C'$ so that
\[ C'_b = \Jac(S_b)[\ell] \setminus \{0 \}. \]
The curve thus redefined has now a good chance of being geometrically connected, and will contain\footnote{Unless~$\rho$ comes from~$\H^1(B_{\overline \Q},\Z/\ell\Z)$, but in this case we could compute it in the Jacobian of~$B$ in view of~\eqref{H_complete}.}~$\rho$ in its~$\H^1$. Assuming that~$C'$ is indeed connected, let~$C$ be the smooth proper model of~$C'$ over~$\Q$; as~$\H^1(C'_{\overline \Q},\Z/\ell \Z)$ is an extension of~$\H^1(C_{\overline \Q},\Z/\ell \Z)$ by copies of~$\mu_\ell^\vee$ according to~\eqref{H_incomplete}, we conclude by~\eqref{H_complete} that the twist of~$\rho$ by a power of the cyclotomic character will be contained in~$\Jac(C)[\ell]$ up to semi-simplification.

\bigskip
This leads to the following plan of attack to compute~$\rho \subset \H^{2} ({S}_{\overline \Q}, \Z/\ell \Z)$:

\begin{enumerate}
\item Compute the Galois representations afforded by the~$\ell$-torsion of the Jacobian of the fibre~$S_b$ of~$\pi$ for various~$b \in B$,
\item Interpolate to glue these data into an explicit model of the cover~$C \longrightarrow B$,
\item Catch a twist of~$\rho$ in the~$\ell$-torsion of the Jacobian of~$C$.
\end{enumerate}

The first and last steps require one to be able to compute explicitly the representations afforded by the torsion of the Jacobian of any curve, which is possible thanks to the method presented in~\cite{Hensel}. In fact, this is the reason why we invented this method in the first place.

The ``interpolation'' part of the second step, as presented here, is quite vague. Fortunately, this will not be a problem for the example that we have in mind, because the fibres~$S_b$ will be elliptic curves. We have not yet studied how to treat the general case. 

\section{Computation of a nice model of~$C$}\label{sect:curvered}

We now apply the method presented in the previous section to the case of the representation~$\rho_{u,3}$ introduced in section~\ref{sect:intro}. According to~\cite[3.10]{Top}, for all~$\ell$, the mod~$\ell$ representation~$\rho_{u,\ell} \otimes \left( \frac{-2}\cdot \right)$ ought to be contained in~$\H^2(S_{\overline \Q},\Z/\ell\Z)$, where~$S$ is the minimal regular model of the projective closure of the surface over~$\Q$ of equation
\begin{equation} z^2=xy(x^2-1)(y^2-1)(x^2-y^2+2xy). \label{eqn:S} \end{equation}
In order to apply the method presented in the previous section to this surface, we need to choose a non-constant map ~$\pi : S \longrightarrow B$, where~$B$ is a curve. We choose~$B = \P^1$ with coordinate~$\lambda$, and~$\pi$ the map sending~$(x,y,z)$ to~$\lambda=y/x$. 

\begin{rk}
Although~$\pi$ is only a rational map from the surface defined by equation~\eqref{eqn:S} to~$\P^1$, it becomes a morphism once we blow up this surface in order to obtain a regular model. In fact,~$\pi$ is the canonical map of~$S$ (cf.~\cite[3.2]{Top}).
\end{rk}

The fibre of~$\pi$ at~$\lambda$ is obtained by setting~$y = \lambda x$ in~\eqref{eqn:S}. The locus of bad fibres~is
\[ Z =\{0,\pm1,1\pm\sqrt2,\infty\} \subset B,\]
and for~$\lambda \not \in Z$, the fibre is an elliptic curve~$E_\lambda$ over~$\Q$. These~$E_\lambda$ define an elliptic curve over~$\Q(\lambda)$ isomorphic to the twist of
\begin{equation} y^2=(x-2\lambda)(x+2\lambda)(x+\lambda^2+1) \label{eq:El} \end{equation}
by~$\lambda (\lambda^2-2\lambda-1)$, and whose~$j$-invariant is~$\displaystyle 2^4 \frac{(\lambda^4+14\lambda^2+1)^3}{\lambda^2(\lambda^2-1)^4}$.

\begin{rk}
Equation~\eqref{eq:El} reveals that~$E_\lambda[2]$ is already defined over~$\Q(\lambda)$. This reflects the fact that the mod~$2$ representation~$\rho_{u,2}$ attached to~$u$ is trivial (up to semi-simplification), as we noted in remark~\ref{rk:rho2triv}.
\end{rk}

The~$\ell$-division polynomial~$\psi_{\ell,\lambda}(x)$ of~$E_\lambda$ is easily computed thanks to~\cite{gp}. By definition, for each~$\lambda$, the Galois action on the roots of~$\psi_{\ell,\lambda}(x)$ describes the Galois action on the~$x$-coordinates of the points of~$E_\lambda[\ell]$. We then compute~$R_{\ell,\lambda}(y)$, the resultant in~$x$ of~$\psi_{\ell,\lambda}(x)$ and of the Weierstrass equation of~$E_\lambda$, which yields a polynomial describing the Galois action on the~$y$-coordinates of~$E_\lambda[\ell]$. For~$\ell=3$, the~$y$-coordinate happens to be injective on~$E_\lambda[\ell]$ for generic~$\lambda$; indeed,~$R_{3,\lambda}(y)$ is squarefree for~$\lambda=2$.

\bigskip

We have thus computed the mod 3 Galois representation afforded by the Jacobian of the fibre of~$\pi$ in terms of~$\lambda$. Substituting~$x$ for~$\lambda$, we obtain the following rather ugly plane model for~$C$:

\vspace{-2mm}

\tiny
\begin{multline*}
-256x^{56} + 6144x^{55} - 62464x^{54} + 333824x^{53} - 859648x^{52} - 120832x^{51} + 7252992x^{50} - 16046080x^{49} - 9891072x^{48} + 90136576x^{47} \\
- 73076736x^{46} - 237805568x^{45} + 420485120x^{44} + 341843968x^{43} - 1165840384x^{42} - 192667648x^{41} + 2178936320x^{40} - 238563328x^{39} \\
- 3063240704x^{38} + 639488000x^{37} + 3412593664x^{36} - 639488000x^{35} - 3063240704x^{34} + 238563328x^{33} + 2178936320x^{32} + 192667648x^{31}\\
 - 1165840384x^{30} - 341843968x^{29} + (-288y^{4} + 420485120)x^{28} + (3456y^{4} + 237805568)x^{27} + (-14400y^{4} - 73076736)x^{26} \\ + (14976y^{4} - 90136576)x^{25} + (56160y^{4} - 9891072)x^{24} + (-142848y^{4} + 16046080)x^{23} + (-52992y^{4} + 7252992)x^{22} + (400896y^{4} + 120832)x^{21} \\
 + (-55872y^{4} - 859648)x^{20} + (-624384y^{4} - 333824)x^{19} + (134784y^{4} - 62464)x^{18} + (624384y^{4} - 6144)x^{17} + (-55872y^{4} - 256)x^{16} \\
 + (16y^{6} - 400896y^{4})x^{15} + (-96y^{6} - 52992y^{4})x^{14} + (-384y^{6} + 142848y^{4})x^{13} + (3232y^{6} + 56160y^{4})x^{12} + (-5424y^{6} - 14976y^{4})x^{11} \\ + (960y^{6} - 14400y^{4})x^{10} - 3456y^{4}x^{9} + (960y^{6} - 288y^{4})x^{8} + 5424y^{6}x^7 + 3232y^6x^6 + 384y^6x^5 - 96y^6x^4 - 16y^6x^3 + 27y^8= 0.
\end{multline*}
\normalsize
A~\cite{Magma} session still manages to reveal in a few seconds that~$C$ is geometrically integral and has (geometric) genus~$g=7$. This is good news, as our method~\cite{Hensel} probably cannot reasonably cope with genera beyond~$20$ or~$30$. 

\begin{rk}\label{rk:Kodaira}
We can easily do the same computation for other values of~$\ell$, and thus get plane models of curves~$C$ that ought to contain a twist of~$\rho_{u,\ell}$ in their Jacobian. However, already for~$\ell=5$, the model we get is so terrible that~\cite{Magma} is unable to determine its genus (the computation was interrupted after 5 days, because it was using more than 400GB of RAM).

We can still compute this genus, by exploiting the fact that~$\pi : S \longrightarrow B$ is an elliptic surface. Indeed, since~$B = \P^1_\lambda$ has genus 0, Riemann-Hurewicz tells us that if~$C$ is connected, then its genus is
\[ g = 1-d+\frac12 \sum_{c \in C} (e_c-1), \]
where~$d=\ell^2-1$ is the degree of the projection~$\psi : C \longrightarrow B$ induced by~$\pi$, and the~$e_c$ are its ramification indices. 

Rewrite
\[ \sum_{c \in C} (e_c-1) = \sum_{\lambda \in B} \sum_{\psi(c) = \lambda} (e_c-1), \]
and notice that for each~$\lambda$, we have
\[ \sum_{\psi(c) = \lambda} (e_c-1) = \sum_{\psi(c) = \lambda} e_c - \sum_{\psi(c) = \lambda} 1 = d - \#\psi^{-1}(\lambda). \]
Besides, the ramification of~$\psi$ can only come from the bad fibres of~$\pi$, so this expression is~$0$ for~$\lambda \not \in Z$.

Our surface~$S$ is the minimal proper regular model of~$E_\lambda / B$, so we can analyse its bad fibres thanks to Tate's algorithm. It reveals that at~$\lambda=0$ and~$\infty$, the special fibre is of Kodaira type~$\text{I}_2^*$, which, according to Table~6 of~\cite[V.10]{ell_monodrom}, implies that the monodromy around~$\lambda$ acts on the homology of the fibre of~$\pi$ by~$T = - \smat{1}{2}{0}{1}$. If~$\ell$ is an odd prime, this means that~$\# \psi^{-1}(\lambda)$, which is the number of orbits of~$T$ acting on~$\Fl^2 \setminus \{0 \}$, is~$\frac1{2\ell} \big( 1 (\ell^2-1) + (\ell-1)(\ell-1) \big) = 2 \ell -2$ by Burnside's formula. Similarly, at~$\lambda=\pm1$, the special fibre is of type~$\text{I}_4$, whence~$T= \smat{1}{4}{0}{1}$ and~$\#\psi^{-1}(\lambda)=2\ell-2$; whereas at~$\lambda=1\pm\sqrt2$, the special fibre is of type~$\text{I}_0^*$, whence~$T=- \smat{1}{0}{0}{1}$ and~$\#\psi^{-1}(\lambda)=\frac{\ell^2-1}2$.

As a result, we find that if~$\ell$ is an odd prime and if~$C$ is connected, then its genus is
\[ g = \frac32 \ell^2-3\ell+\frac52. \]
In particular, we recover~$g=7$ for~$\ell=3$, and we find that~$g=25$ for~$\ell=5$ and that~$g=55$ for~$\ell=7$. This means that our method~\cite{Hensel} could probably manage to compute the mod 5 representation~$\rho_{u,5}$ if we could find a decent enough model for~$C$ for~$\ell=5$ and if we were patient enough, whereas~$\ell \geqslant 7$ seems out of our reach.
\end{rk}

Let us get back to the case~$\ell=3$ and to our curve~$C$ of genus~$7$. The model that we have just obtained has degree 56, and therefore arithmetic genus 1485. We do not want to work with such a badly singular model, so we attempt to eliminate the worst of the singularities by having~\cite{Magma} determine the canonical image of~$C$ in~$\P^6$ and project it on a plane. This yields the already more appealing model

\vspace{-2mm}

\tiny
\begin{multline*}
(-374594220y^6 + 148459311y^5 - 20961720y^4 + 1285362y^3 - 35100y^2 + 351y)x^4 \\ + (-61958809438y^8 + 12030741624y^7 - 574743724y^6 - 5928484y^5 + 27600y^4 + 129884y^3 - 8516y^2 + 216y - 2)x^2 \\
+ (15790199962940y^{10} - 5854413418867y^9 + 927447207596y^8 - 81010188948y^7 + 4049824636y^6- 100135334y^5\\ - 48724y^4 + 70252y^3 - 1664y^2 + 13y) = 0.
\end{multline*}
\normalsize
We check that this model still has genus 7, which by Riemann-Hurwitz ensures that it is birational to the previous one, as opposed to being a quotient of it.

\bigskip

By projecting the canonical image onto another plane, we also find that~$C$ is a cover of degree 2 (simply given by~$x \mapsto x^2$ on our new model) of a curve~$E$ of genus~1, which turns out to be an elliptic curve isomorphic to the modular curve~$X_0(24)$. We learn from the~\cite{LMFDB} that this curve has rank 0; after careful back-tracking, this allows us to determine the complete list of rational points of~$C$ in our new model.

\bigskip

In order to further simplify our model for~$C$, we turn to a trial-and-error manual shifting and rescaling process based on the shape of the rational points thus obtained and on the observation of the~$p$-adic valuations of the coefficients of our model for small primes~$p$. After a few attempts, we obtain the quite satisfying model
\begin{multline*}
(3y^5 - 6y^3 + 3y)x^4 + (2y^8 - 8y^7 + 4y^6 + 12y^5 + 12y^3 - 4y^2 - 8y - 2)x^2\\ + (9y^9 - 36y^8 - 36y^7 + 36y^6 + 18y^5 - 36y^4 - 36y^3 + 36y^2 + 9y) = 0.
\end{multline*}

\begin{rk}\label{rk:curvered}
The whole process to obtain a nice model for~$C$ was rather rustic. Most computer algebra systems include algorithms that, given a complicated polynomial defining a number field, are able to find a much simpler polynomial defining the same field (when it exists); it would be nice to have similar algorithms for curves!
\end{rk}

The arguments of the previous section show that the 3-torsion of the Jacobian of this curve contains the representation contained in~$\H^2(S_{\overline \Q},\Z/3\Z)$ up to twist by the mod~$3$ cyclotomic character~$\chi_3$, which agrees with the quadratic character~$\left( \frac{-3} \cdot \right)$. Since this representation was already the twist of the representation~$\rho_{u,3}$ we are interested in by~$\left( \frac{-2} \cdot \right)$, this is just an extra twist. 

\bigskip

In order to confirm this, we can check that the characteristic polynomials match at a few primes~$p$. Indeed, let~$\rho'_{u,3}$ be the~$\GL_6(\F_3)$-valued representation obtained by restricting the scalars from the~$\GL_3(\F_9)$-valued representation~$\rho_{u,3}$. On the one hand, we know that  for~$p \neq 2,3$, the characteristic polynomial of~$\rho'_{u,3}(\Frob_p)$ is
\[ \chi'_p(x) = \chi_p(x) \overline{\chi_p(x)} = (x^3-a_p x^2+p \overline{a_p} x-p^3)(x^3-\overline{a_p} x^2+p a_p x - p^3) \in \F_3[x], \]
the norm from~$\F_9$ to~$\F_3$ of the polynomial~$\chi_p(x)$ given in~\eqref{eqn:chip}, and furthermore~\cite[2.5]{Top} provides us with the values of the Hecke eigenvalues~$a_p$ of~$u$ for~$p \le 67$. On the other hand, we can determine the characteristic polynomial~$L_p(x) \in \Z[x]$ of the Frobenius at~$p$ acting on the Jacobian of~$C$ (which is the local factor at~$p$ of its L~function) by counting the~$\F_{p^a}$-points of~$C$ for~$a \le g$, where~$g = 7$ is the genus of~$C$; in practice,~\cite{Magma} can do this in reasonable time for~$p \le 19$. We then check that for~$5 \le p \le 19$,~$L_p(x) \bmod 3$ is divisible by the characteristic polynomial~$\chi'_p(\epsilon x)$ of the image of~$\Frob_p$ by~$\rho'_{u,3} \otimes \left( \frac{6} \cdot \right)$, where~$\epsilon = \big( \frac{6} p \big) = \pm 1$. This confirms our hopes that the Jacobian~$J$ of~$C$ contains~$\rho'_{u,3} \otimes \left( \frac{6} \cdot \right)$ in its 3-torsion.

\begin{rk}
Since~$g=7$, the degree of~$L_p(x)$ is 14. We actually observe that for all the primes~$5 \le p \le 19$ that we can test,~$L_p(x)$ is congruent mod 3 to the product
\[ \chi_{E,p}(x) \, \chi'_p\big( (\scalebox{0.8}{$\frac{6}p$}) x \big) \, \chi'_p\big( (\scalebox{0.8}{$\frac{-2}p$}) x \big) \]
where the factors have respective degrees 2, 6, and 6, and are the characteristic polynomial of~$\Frob_p$ for the mod~$3$ representation~$\rho_{E,3}$ attached to the elliptic \linebreak curve~$E=X_0(24)$ exhibited above, the expected twist~$\rho'_{u,3} \otimes \left( \frac{6} \cdot \right)$, and the twist~$\rho'_{u,3} \otimes \left( \frac{-2} \cdot \right)$ originally contained in~$\H^2(S_{\overline \Q},\Z/3\Z)$, respectively. The fact that we eventually managed to compute~$\rho'_{u,3} \otimes \left( \frac{6} \cdot \right)$ from a piece of~$J[3]$, whereas we were unable to do the same for~$\rho'_{u,3} \otimes \left( \frac{-2} \cdot \right)$ (even after significantly increasing the~$p$-adic accuracy in our computation, cf. the next section), leads us to guess that the Galois-module~$J[3]$ decomposes as
\[ J[3] \sim \left[ \begin{matrix} \rho_{E,3} & & \\ & \rho'_{u,3} \otimes \left( \frac{6} \cdot \right) & * \\ & & \rho'_{u,3} \otimes \left( \frac{-2} \cdot \right) \end{matrix} \right], \]
where~$*$ is non-trivial. In other words, the twist~$\rho'_{u,3} \otimes \left( \frac{-2} \cdot \right)$ originally contained in~$\H^2(S_{\overline \Q},\Z/3\Z)$ seems to show up as \emph{quotient} of~$J[3]$. 
\end{rk}

\section{Computation of the representation in~$J$}\label{sect:Hensel}

Now that we have obtained a reasonable model for~$C$, we are going to use our method~\cite{Hensel} to compute the representation~$\rho'_{u,3} \otimes \left( \frac{6} \cdot \right)$ afforded by a Galois submodule~$T$ of the 3-torsion of the Jacobian~$J$ of~$C$. This method requires us to pick a prime~$p$ of good reduction for which the local L factor
\[ L_p(x) = \det(x-\Frob_p {\vert J}) \in \Z[x]\]
and the characteristic polynomial
\[\chi'_p\big( (\scalebox{0.8}{$\frac{6}p$}) x \big) = \det(x - \Frob_p \vert T) \in \F_3[x]\]
are known, and such that~$\chi'_p\big( (\scalebox{0.8}{$\frac{6}p$}) x \big)$ is coprime to its cofactor~$L_p(x) / \chi'_p\big( (\scalebox{0.8}{$\frac{6}p$}) x \big)$. We choose~$p=11$, which satisfies these assumptions. Besides, we then have
\[ \chi'_p\big( (\scalebox{0.8}{$\frac{6}p$}) x \big) = \sum_{k=0}^6 x^k \in \F_3[x], \]
which is irreducible; this fact will be useful on two occasions below.

Our method performs by generating points of~$J[3]$ over an appropriate extension of~$\F_{11}$ and projecting them onto~$T$ thanks to the action of~$\Frob_{11}$ until we get a basis of~$T$, to lift this basis 11-adically, and then to evaluate a Galois-equivariant rational map~$\alpha$ from~$J$ to~$\A^1$ at the points of~$T$. If the 11-adic accuracy is sufficient, then we will be able to identify the coefficients of the polynomial
\[ F(x) = \prod_{\substack{t \in T \\ t \neq 0}} \big( x - \alpha(t) \big) \]
as rational numbers; if furthermore~$\alpha$ is injective on~$T$, we will thus have obtained a polynomial whose roots are permuted under Galois just as the non-zero points of~$T$ are.

Here, the fact that the characteristic polynomial of~$\Frob_{11}$ is squarefree implies that~$\Frob_{11}$ is a cyclic endomorphism of the~$\F_3$-vector space~$T$, which means that we can afford to lift only one point of~$T$ and recover a basis by applying~$\Frob_{11}$ repeatedly (cf.\ section 6.4 of~\cite{Hensel}).

\begin{rk}\label{rk:bestappr}
The fact that we identify the coefficients of~$F(x)$ from their~$p$-adic approximations is the reason why we cannot rigorously certify that our computation results are correct. Nevertheless, we will give in section~\ref{sect:Dok} below very strong evidence that these results are correct beyond reasonable doubt.
\end{rk} 

\pagebreak

In order to be able to compute in~$J$, we need to fix an effective divisor~$D_0$ on~$C$ of \linebreak degree~$d_0 \ge 2g+1=15$ for which we can explicitly compute the corresponding Riemann-Roch space. In order to construct the evaluation map~$\alpha$, we also need to pick two non-equivalent effective divisors~$E_1$,~$E_2$ of degree~$d_0-g$ such that we can also compute the Riemann-Roch spaces attached to~$2D_0-E_1$ and~$2D_0-E_2$ (the notations are the same as in~\cite{Hensel}).

It is qualitatively clear that we should strive to choose~$D_0$,~$E_1$, and~$E_2$ so that these three Riemann-Roch spaces are as ``nice'' as possible, since the values of~$\alpha$ will then have smaller arithmetic height, so that the~$p$-adic accuracy required to identify the coefficients of~$F(x)$ will be lower and the computation will be more efficient. After a bit of experimentation with~\cite{Magma}, we choose 
\[ d_0=16, \quad D_0 = 9P+7Q, \quad E_1 = 6P+3Q, \quad E_2=5P+4Q, \]
where~$P,Q \in C(\Q)$ are points such that, in the model obtained at the end of the previous section, the divisors of poles of~$x$ and~$y$ are respectively
\[ (x)_\infty=3P+Q+R+M_1+M_2 \text{ and } (y)_\infty=2P+2Q \]
where~$R$ has degree 1 and~$M_1$ and~$M_2$ both have degree 2.

\begin{rk}
It may happen that the~$\Q$-basis of a Riemann-Roch space becomes linearly dependent when reduced mod~$p$. Fortunately, this is easy to detect, because functions on~$C$ are represented internally in~\cite{Hensel} as the vector of their values at a large enough set of fixed points of~$C$. This is also easy to fix, by Gaussian elimination: given~$s_1, \cdots, s_d \in \Q(C)$ forming the basis of a given Riemann-Roch space, if~$\sum_i \lambda_i s_i \equiv 0 \bmod p$ for some~$\lambda_i \in \F_p$ not all 0, it suffices to substitute~$\frac1p \sum_i \widetilde{\lambda_i} s_i$ to~$s_j$, where~$j$ is such that~$\lambda_j \neq 0$ and the~$\widetilde \lambda_i$ are lifts to~$\Z$ of the~$\lambda_i$. However, this complicates the basis of the Riemann-Roch space, which in turn increases the height of the values of the evaluation map~$\alpha$. Fortunately, this phenomenon does not occur with our choices of~$D_0$,~$E_1$,~$E_2$ and~$p$.
\end{rk}

\bigskip

Now that we have made these choices, we are ready to launch the computations. After about 30 hours of CPU time (but only 1 hour of real time, thanks to parallelisation), we obtain a polynomial~$F(x)$ of degree~$3^6-1=728$ whose coefficients are rational numbers which all have (up to some small factors) the same denominator, a 191-digit integer. The~$p$-adic precision used was~$O(11^{1024})$.

The discriminant of~$F(x)$ factors into a large power of 2 times a huge power of 3 times a large square, which indicates that its coefficients have probably been correctly identified from their~$p$-adic approximations. The fact that discriminant is non-zero also shows that~$\alpha$ is injective on~$T$ minus the origin, so that the roots of~$F(x)$ represent faithfully the Galois action on~$T$ minus the origin, as desired.

\section{The image of~$\rho_{u,3}$}\label{sect:image}

We find that the polynomial~$F(x)$ computed in the previous section has three factors over~$\Q$, of respective degrees 224, 252, and 252. This shows that the image of our representation does not act transitively on~$\F_3^6$ minus the origin.

However, the degrees of these factors do not clearly indicate which subgroup of~$\GL_6(\F_3)$ we are dealing with. In order to figure this out, we would like to determine the Galois groups of these factors; however, their degrees and heights are far too large for standard Galois group computation algorithms. We would therefore like to reduce these factors (in the sense of remark~\ref{rk:curvered}), but they are actually too large even for this!

\bigskip

As in~\cite[section 2]{certif}, we circumvent this problem by considering the \emph{projective} version of our representation, which has values in~$\PGL_3(\F_9)$. We can obtain a polynomial corresponding to this representation by gathering the 11-adic roots~$\alpha(t)$ of~$F(x)$ along the~$\F_9$-vector lines of~$T$ in a symmetric way (e.g. by summing or multiplying them). However, this requires us to understand the~$\F_9$-structure of~$T$, whereas we only know the~$\F_3$-structure for now.

Let~$\Phi \in \GL(T)$ be the action of~$\Frob_{11}$ on~$T$. We know that~$\F_9^\times$ acts on~$T$ by a cyclic subgroup of~$\GL(T)$ of order 8 contained in the commutant of~$\Phi$.  Luckily,~$\Phi$ is cyclic, so its commutant is simply~$\F_3[\Phi]$, which is a ring isomorphic to~$\F_{3^6}$ since the characteristic polynomial of~$\Phi$ is isomorphic over~$\F_3$. In particular, there is a unique cyclic subgroup of order~8 in~$\F_3[\Phi]^\times$.

We can thus compute as above a polynomial~$F_0(x)$ of degree~$\frac18 \deg F = 91$ describing the projective representation attached to~$\rho_{u,3} \otimes \left( \frac{6} \cdot \right)$ (which is also that attached to~$\rho_{u,3}$).

\bigskip

We can also be a bit more subtle, and consider all intermediate representations between the linear one and the projective one. Let us write~$\F_9 = \F_3(i)$, where~$i^2=-1$. Then the subgroups of~$\F_9^\times$ are, in decreasing order, 
\[ \F_9^* \ \ge \ \{\pm1, \pm i\} \ \ge \ \{\pm1\} \ \ge \ \{1\}.\]
We can construct as above polynomials~$F_0(x)$,~$F_1(x)$,~$F_2(x)$ and~$F_3(x)=F(x)$ describing the corresponding quotients of~$\rho_{u,3} \otimes \left( \frac{6} \cdot \right)$. These polynomials factor over~$\Q$ as follows:
\[ \begin{array}{c|c} 
\text{ Quotient by } & \text{ Degrees of factors } \\
\hline
\F_9^\times & 28 + 63 \\
\{ \pm1,\pm i\} & 56+63+63 \\
\{ \pm 1 \} & 112+126+126 \\
\{1\} & \phantom{.}224+252+252. \\
\end{array} \]
These degrees still do not clearly indicate what the image of the representation is. If anything, the fact that we have two factors for the projective representation that become three factors afterwards is rather mysterious.

\bigskip

Fortunately, one of the factors of~$F_0(x)$ has degree 28, which is small enough that we can compute a much simpler polynomial defining the same number field, namely
\begin{multline*}
x^{28} - 12x^{27} + 60x^{26} - 132x^{25} - 30x^{24} + 624x^{23} + 420x^{22} - 7704x^{21} \\ + 17118x^{20} - 9504x^{19} - 14424x^{18} + 10824x^{17} + 36492x^{16} - 64992x^{15} + 19488x^{14} \\ + 56064x^{13} - 89604x^{12} + 109296x^{11} - 88368x^{10} - 11472x^9 + 58488x^8 - 130176x^7 \\ + 34224x^6 - 58272x^5 - 39960x^4 + 32256x^3 + 24480x^2 - 352x - 1776.
\end{multline*}
This polynomial is nice enough that~\cite{Magma} can rigorously determine its Galois group in less than a minute. This group turns out to be the simple group~$\PSUt$, which explains all the observations made above!

\bigskip

Indeed, first of all one checks thanks to~\cite{gp} that the field defined by the factor of degree 63 of~$F_0(x)$ is contained in the compositum of the field defined by that of degree 28 with itself, which shows that these factors have the same splitting field. Next, since there are no nontrivial cube roots of unity in characteristic 3, the quotient~$\SUt \longrightarrow \PSUt$ is actually an isomorphism; in particular, it admits a section. This means that one twist of~$\rho_{u,3}$ has image~$\SUt$ (and actually, this twist is the one by the mod 3 cyclotomic character~$\chi_3 = \left( \frac{-3} \cdot \right)$ since we have seen that~$\det \rho_{u,3} = \chi_3^3 = \chi_3$).

Let now~$H$ be a non-degenerate Hermitian form on the space~$\F_{q^2}^n$, where~$q$ is a prime power and~$n \in \N$, and let~$A_n$ (respectively~$B_n$) be the number of elements~$t \in \F_{q^2}^n$ such that~$H(t)=1$ (respectively, such that~$H(t)=0$). Since the norm between finite fields is surjective,~$A_n$ is also the number of elements~$t \in \F_{q^2}^n$ such that~$H(t)$ has prescribed value~$y \in \F_q^\times$. From this fact, one easily determines a crossed recurrence relation satisfied by~$A_n$ and~$B_n$, from which one deduces that
\[ A_n = q^{2n-1}+(-q)^{n-1}, \quad B_n = q^{2n-1} - (q-1)(-q)^{n-1}. \]
For~$n=3$ and~$q=3$, one finds~$A_n=252$ and~$B_n=225$, which explains the shape~$224+252+252$ of the factorization of~$F_3(x)$: the first factor corresponds to the nonzero isotropic~$t \in T$, and the other two correspond to the~$t$ such that~$H(t)=1$ (respectively, such that~$H(t)=-1$).

Similarly, the Galois group of the factor of degree 28 of~$F_0(x)$ is per\-mu\-ta\-tion-iso\-mor\-phic to~$\PSUt$ acting on the isotropic lines of~$\F_9^3$, whereas that of the factor of degree 63 corresponds to the action of~$\PSUt$ on non-isotropic lines.

Finally, the fact that the value of~$H$ is not well defined at a non-isotropic~$t$ known up to scaling by~$\F_9^\times$, but becomes well-defined if we know~$t$ up to scaling by
\[ \{\pm1, \pm i \} = \Ker \operatorname{Norm} : \F_9^\times \rightarrow \F_3^\times, \]
explains why the factor of degree 63 of~$F_0(x)$ yields two factors of degree 63 of~$F_1(x)$ instead of one of degree 126.

\begin{rk}
We have also obtained a simpler polynomial of degree 63 defining the same number field as this factor of degree 63. This polynomial is available on the author's web page~\cite{webpage}, and we do not reproduce it here. The polynomial of degree 28 displayed above and this polynomial of degree 63 thus solve the inverse Galois problem for the standard actions of the simple group~$\PSUt$ in respective degrees 28 and 63, and with controlled ramification (only at 2 and 3) to boost! The respective discriminants and signatures of the corresponding number fields are as follows:
\[ \begin{array}{c|c|c} 
\text{Degree} & \text{Discriminant} & \text{Signature} \\
\hline
\phantom{\Big\vert} 28 \phantom{\Big\vert} & 2^{76} 3^{48} & (4,12) \\ 
63 & 2^{166} 3^{108} & \phantom{.}(7,28). \\ 
\end{array} \]
\end{rk}

\bigskip

\begin{rk}
Since factors of the~$F_i(x), \ 0 \le i \le 3$ correspond to towers of quadratic extensions, we can use the techniques presented in~\cite[section 2]{certif} to compute nice polynomials defining the same number fields as the factors of~$F_3(x)$. This technique has the advantage of naturally producing even polynomials, such that the field automorphism induced by~$x \mapsto -x$ corresponds to the action of~$-1 \in \F_9^\times$. This means that given such a polynomial~$f(x^2)$, the polynomial~$f(Dx^2)$ corresponds to the twist of the representation by~$\left( \frac{D} \cdot \right)$ for any~$D \in \Q^\times$. By taking~$D=-3$, we obtain polynomials corresponding to the representation~$\rho_{u,3} \otimes \left( \frac{-3} \cdot \right)$ whose image is the simple group~$\SUt \simeq \PSUt$, thus again solving the inverse Galois problem for the natural permutation representations of this group. These polynomials are also available on the author's web page~\cite{webpage}.
\end{rk}

\pagebreak

\section{Computation of the image of~$\Frob_p$}\label{sect:Dok}

Let~$L$ be a Galois number field, given as the splitting field of an irreducible \linebreak polynomial~$f(x) \in \Q[x]$. In~\cite{Dok}, the Dokchitsers show that if the action of~$G = \Gal(L/\Q)$ on the roots of~$f(x)$ in~$L$ is known explicitly, then one can compute pairwise coprime resolvents~$\Gamma_C(x) \in \Q[x]$ indexed by the conjugacy classes~$C$ of~$G$, such that if the image of~$\Frob_p$ in~$G$ lies in~$C$, then the corresponding resolvent~$\Gamma_C(x)$ vanishes at~$x_p=\operatorname{Tr}^{A_p}_{\F_p} \big( a^p h(a) \big) \in \F_p$, where~$A_p = \F_p[x]/f(x)$,~$a$ is the image of~$x$ in~$A_p$, and~$h(x) \in \Z[x]$ is a fixed parameter on which the~$\Gamma_C(x)$ depend.

The point is that since the~$\Gamma_C(x)$ are coprime, they remain coprime mod~$p$ for almost all~$p$, so only one of them can vanish at~$x_p$ and we can tell in which class~$C$ the image of~$\Frob_p$ lies. The finitely many~$p$ for which this is no longer true are usually quite small, and for these~$p$ be get not one but several~$C$ that may contain~$\Frob_p$. If the conjugacy class of~$\Frob_p$ is really wanted for such a~$p$, one should recompute the resolvents with another value of the parameter~$h$.

\bigskip

Since the quotient~$\SUt \longrightarrow \PSUt$ is actually an isomorphism, and as~$\det \rho_{u,3} = \left( \frac{-3} \cdot \right)$ is known explicitly, for each prime~$p$ we can recover the image of~$\Frob_p$ by~$\rho_{u,3}$ from its image by the projective version of this representation. Namely, if the image of~$\Frob_p$ by the projective representation is conjugate to~$\overline M \in \PSUt$, then~$\rho_{u,3}(\Frob_p)$ is conjugate to~$\big( \frac{-3} p \big) M$ in~$\Ut$, \linebreak where~$M \in \SUt$ is the image of~$\overline M$ by the inverse of this isomorphism.

\bigskip

We thus apply the Dokchitsers' method to the case where~$f(x)$ is the polynomial of degree~$28$ displayed in the previous section. This polynomial has Galois group~$\PSUt$, and its roots are indexed by the lines of~$\F_9^3$ that are isotropic with respect to a certain hermitian form~$H$. We can determine~$H$ from the fact that it is preserved by the action of~$\Frob_{11}$, and after a change of basis of~$\F_9^3$ we can assume that~$H$ is the standard Hermitian form.

The group~$\PSUt$ has order 6048, which is small enough that~\cite{Magma} can effortlessly decompose it explicitly into conjugacy classes, which is all we need to compute the resolvents~$\Gamma_C(x)$.

\bigskip

Thanks to these resolvents, we can now determine the image of~$\Frob_p$ by~$\rho_{u,3}$ for almost all~$p$. Let us start by the primes between 5 and 67, for which the value of the Hecke eigenvalue~$a_p \in \Z[i]$ is given in~\cite{Top}. 

\newpage

\rowcolors{2}{gray!25}{white}
\[ 
\begin{array}{|c!{\rule[-3ex]{0ex}{7ex}}cc|}
\hline
\rowcolor{gray!50}
\phantom{abc} p \phantom{abc} & \phantom{abcdefghi} \rho_{u,3}(\operatorname{Frob}_p)  \phantom{abc} & \phantom{ab} a_p \text{ from~\cite{Top}} \phantom{ab} \\
5 & \phantom{abcdei} 3 \text{ possibilities}  & -4i-1 \\
7 &  \phantom{abcd} +\left[ \begin{Smallmatrix} 0 & i + 1 & i - 1 \\ 0 & i + 1 & -i + 1 \\ 1 & 0 & 0 \end{Smallmatrix} \right] & 4i+1 \\
11 &  \phantom{abcd} -\left[ \begin{Smallmatrix} 0 & i + 1 & i - 1 \\ 0 & i + 1 & -i + 1 \\ 1 & 0 & 0 \end{Smallmatrix} \right] & -10i-7 \\
13 &  \phantom{abcd} +\left[ \begin{Smallmatrix} 0 & i + 1 & i + 1 \\ 0 & i - 1 & -i + 1 \\ 1 & 0 & 0 \end{Smallmatrix} \right]  & 4i-1 \\
17 &  \phantom{abcd} -\left[ \begin{Smallmatrix} 1 & 0 & 0 \\ 0 & i - 1 & i - 1 \\ 0 & i +1 & -i - 1 \end{Smallmatrix} \right] & 7 \\
19 &  \phantom{abcd}+\left[ \begin{Smallmatrix} 0 & i + 1 & i - 1 \\ 0 & i + 1 & -i + 1 \\ 1 & 0 & 0 \end{Smallmatrix} \right] & -14i+1 \\
23 &  \phantom{abcd} -\left[ \begin{Smallmatrix} 0 & i + 1 & i - 1 \\ 0 & i + 1 & -i + 1 \\ 1 & 0 & 0 \end{Smallmatrix} \right] & -4i+17 \\
29 & \phantom{abcd}  -\left[ \begin{Smallmatrix} 0 & 0 &1 \\ 1 & 0 & 0 \\ 0 & 1 & 0 \end{Smallmatrix} \right]  & -12i-9 \\
31 & \phantom{abcdei} 3 \text{ possibilities}  & 1 \\
37 & \phantom{abcdei} 2 \text{ possibilities}  & 28i-25 \\
41 & \phantom{abcd} -\left[ \begin{Smallmatrix} 1 & 0 & 0 \\ 0 & i - 1 & i - 1 \\ 0 & i + 1 & -i - 1 \end{Smallmatrix} \right]  & -5 \\
43 & \phantom{abcd} +\left[ \begin{Smallmatrix} 1 & 0 & 0 \\ 0 & i - 1 & i - 1 \\ 0 & i + 1 & -i - 1 \end{Smallmatrix} \right]  & 30i-7 \\
47 & \phantom{abcd} -\left[ \begin{Smallmatrix} 0 & i + 1 & -i - 1 \\ 0 & -i + 1 & -i + 1 \\ 1 & 0 & 0 \end{Smallmatrix} \right]  & 40i+17 \\
53 &  \phantom{abcdei} 2 \text{ possibilities}  & -20i+23 \\
59 & \phantom{abcd} -\left[ \begin{Smallmatrix} 0 & 0 & -i \\ 0 & -i & 0 \\ 1 & 0 & 0 \end{Smallmatrix} \right]  & 22i-39 \\
61 & \phantom{abcd} +\left[ \begin{Smallmatrix} 0 & 0 & -i \\ 0 & -i & 0 \\ 1 & 0 & 0 \end{Smallmatrix} \right]  & 20i+63 \\
67 & \phantom{abcd} +\left[ \begin{Smallmatrix} 0 & i + 1 & -i + 1 \\ 0 & -i - 1 & -i + 1 \\ 1 & 0 & 0 \end{Smallmatrix} \right]  & -22i+65 \\
\hline
\end{array}
\]

\newpage

With our choice our~$h$, the resolvents~$\Gamma_C(x)$ do not remain coprime mod~$p$ \linebreak for~$p \in \{ 5, 31, 37,53 \}$, so for these~$p$ we cannot determine the image of~$\Frob_p$ without recomputing the resolvents with another value of~$h$. For the other~$p$, we can compute the conjugacy class~$C \subset \PSUt$ containing~$\Frob_p$, and we display the image of~$\Frob_p$ by~$\rho_{u,3}$ as
\[ \left( \frac{-3} p \right) M \in \Ut,\]
where~$M$ is a fixed representative of its conjugacy class in~$\SUt \simeq \PSUt$ that we have arbitrarily chosen because many of its coefficients were 0. 

The fact that the trace agrees with the reduction mod 3 of the value of~$a_p$ given in~\cite{Top} is convincing evidence that we have correctly computed~$\rho_{u,3}$.

\bigskip

Next, we do the same thing for the first twenty primes above~$10^{1000}$. Of course, the~$\Gamma_C(x)$ \linebreak remain coprime mod~$p$ for such large~$p$, so we find unambiguously the conjugacy class \linebreak of~$\rho_{u,3}(\Frob_p)$. By looking at the trace, we deduce the value of~$a_p$ mod 3. The results are displayed in the table below.

\begin{rk} 
It takes about 100 seconds for~\cite{gp} to certify the primality of such a large prime, but only 4 seconds to compute the conjugacy class of~$\rho_{u,3}(\Frob_p)$, almost all of which are spent computing~$\operatorname{Tr}^{A_p}_{\F_p} \big( a^p h(a) \big) \in \F_p$.
\end{rk}

\bigskip

The resolvents~$\Gamma_C(x)$ are available on the author's web page~\cite{webpage}.

\newpage

\rowcolors{2}{gray!25}{white}
\[ 
\begin{array}{|c!{\rule[-2.7ex]{0ex}{6.6ex}}cc|}
\hline
\rowcolor{gray!50}
\phantom{abcdef} p \phantom{abcdef} & \phantom{abcdefgh}  \rho_{u,3}(\operatorname{Frob}_p)  \phantom{abcdef} & a_p \bmod 3 \Z[i] \\
10^{1000}+453 &  +\left[ \begin{Smallmatrix} 1 & 0 & 0 \\ 0 & i-1 & i-1 \\ 0 & i+1 & -i-1 \end{Smallmatrix} \right]  & -1 \\
10^{1000}+1357 &  -\left[ \begin{Smallmatrix} 0 & 0 & i \\ 0 & i & 0 \\ 1 & 0 & 0 \end{Smallmatrix} \right]  & -i \\
10^{1000}+2713 &  -\left[ \begin{Smallmatrix} 0 & 0 & -i \\ 0 & -i & 0 \\ 1 & 0 & 0 \end{Smallmatrix} \right]  & i \\
10^{1000}+4351 &  -\left[ \begin{Smallmatrix} 0 & i+1 & -i-1 \\ 0 & -i+1 & -i+1 \\ 1 & 0 & 0 \end{Smallmatrix} \right]  & i-1 \\
10^{1000}+5733 &  +\left[ \begin{Smallmatrix} 0 & i + 1 & -i + 1 \\ 0 & -i - 1 & -i + 1 \\ 1 & 0 & 0 \end{Smallmatrix} \right]  & -i - 1 \\
10^{1000}+7383 &  +\left[ \begin{Smallmatrix} 0 & 0 & -i \\ 0 & -i & 0 \\ 1 & 0 & 0 \end{Smallmatrix} \right]  & -i \\
10^{1000}+10401 & +\left[ \begin{Smallmatrix} 1 & 0 & 0 \\ 0 & i & 0 \\ 0 & 0 & -i \end{Smallmatrix} \right]  & 1 \\
10^{1000}+11979 &  +\left[ \begin{Smallmatrix} 0 & i + 1 & i + 1 \\ 0 & i - 1 & -i + 1 \\ 1 & 0 & 0 \end{Smallmatrix} \right]  & i - 1 \\
10^{1000}+17557 &  -\left[ \begin{Smallmatrix} 1 & 0 & 0 \\ 0 & i - 1 & i - 1 \\ 0 & i + 1 & -i - 1 \end{Smallmatrix} \right]  & 1 \\
10^{1000}+21567 &  +\left[ \begin{Smallmatrix} 1 & 0 & 0 \\ 0 & i - 1 & i - 1 \\ 0 & i + 1 & -i - 1 \end{Smallmatrix} \right]  & -1 \\
10^{1000}+22273 &  -\left[ \begin{Smallmatrix} 0 & i + 1 & -i - 1 \\ 0 & -i + 1 & -i + 1 \\ 1 & 0 & 0 \end{Smallmatrix} \right]  & i - 1 \\
10^{1000}+24493 &  -\left[ \begin{Smallmatrix} 0 & i+ 1 & -i - 1 \\ 0 & -i + 1 & -i + 1 \\ 1 & 0 & 0 \end{Smallmatrix} \right]  & i - 1 \\
10^{1000}+25947 &  +\left[ \begin{Smallmatrix} 0 & i + 1 & i - 1 \\ 0 & i + 1 & -i + 1 \\ 1 & 0 & 0 \end{Smallmatrix} \right]  & i + 1 \\
10^{1000}+27057 &  +\left[ \begin{Smallmatrix} 0 & i + 1 & -i + 1 \\ 0 & -i - 1 & -i + 1 \\ 1 & 0 & 0 \end{Smallmatrix} \right]  & -i - 1 \\
10^{1000}+29737 &  -\left[ \begin{Smallmatrix} 0 & i + 1 & -i + 1 \\ 0 & -i - 1 & -i + 1 \\ 1 & 0 & 0 \end{Smallmatrix} \right]  & i + 1 \\
10^{1000}+41599 &  -\left[ \begin{Smallmatrix} 1 & 0 & 0 \\ 0 & i & 0 \\ 0 & 0 & -i \end{Smallmatrix} \right]  & -1 \\
10^{1000}+43789 &  -\left[ \begin{Smallmatrix} 0 & 0 & 1 \\ 1 & 0 & 0 \\ 0 & 1 & 0 \end{Smallmatrix} \right]  & 0 \\
10^{1000}+46227 &  +\left[ \begin{Smallmatrix} 0 & i + 1 & -i - 1 \\ 0 & -i + 1 & -i + 1 \\ 1 & 0 & 0 \end{Smallmatrix} \right]  & -i + 1 \\
10^{1000}+46339 &  -\left[ \begin{Smallmatrix} 0 & 0 & i \\ 0 & i & 0 \\ 1 & 0 & 0 \end{Smallmatrix} \right]  & -i \\
10^{1000}+52423 &  -\left[ \begin{Smallmatrix} 0 & i + 1 & i - 1 \\ 0 & i + 1 & -i + 1 \\ 1 & 0 & 0 \end{Smallmatrix} \right]  & -i - 1 \\
\hline
\end{array}
\]

\end{document}